\documentclass[11pt]{article} 
\usepackage{amsfonts} 
\usepackage{latexsym} 
\usepackage{graphics} 
\usepackage{psfrag} 
\usepackage{theorem}

\def\bb{\bf}

\theorembodyfont{\rm} 
\newtheorem{theorem}{Theorem}[section] 
\newtheorem{proposition}[theorem]{Proposition} 
\newtheorem{lemma}[theorem]{Lemma} 
\newtheorem{corollary}[theorem]{Corollary} 
\newtheorem{definition}[theorem]{Definition} 
 
\newtheorem{example}[theorem]{Example}

\newcommand{\remembertheorem}[1]{\setcounter{#1}{\value{theorem}}} 
\newcommand{\recalltheorem}[1]{\arabic{section}.\arabic{#1}}

\newcommand{\jremth}[1]{\newcounter{#1theorem} 
\setcounter{#1theorem}{\value{theorem}} 
\newcounter{#1sect} 
\setcounter{#1sect}{\value{section}}} 
\newcommand{\jrecth}[1]{\arabic{#1sect}.\arabic{#1theorem}} 

\newcommand{\Proof}{{\bf Proof: }} 
\newcommand{\QED}{\hfill$\Box$\vskip0.2cm} 
\newcommand{\Complex}{{\hbox{\bb C}}} 
 
\newcommand{\Natural}{{\bb N}}

\newcommand{\stddiff}[1]{D^{1}\left(#1\right)} 
 
\newcommand{\diff}[2]{{\mathcal{D}_{#1}^{1}}\left(#2\right)} 
\newcommand{\diffn}[3]{{\mathcal{D}_{#1}^{#2}}\left(#3\right)} 
\newcommand{\wrt}[1]{\mbox{ }d #1} 
\newcommand{\paths}{\mathcal{F}} 
\newcommand{\Affines}[1]{\mathcal{A}(X)} 
\newcommand{\Smooths}[1]{\mathcal{S}(X)} 
\newcommand{\Rects}[1]{\mathcal{R}(X)} 
 
\newcommand{\Image}[1]{\mbox{Image}\left(#1\right)} 
 
\newcommand{\Closure}[1]{\overline{#1}} 
\newcommand{\Boundary}[1]{\partial #1} 
\newcommand{\Length}[1]{|#1|}

\newcommand{\dis}{\displaystyle}

\def\rd{{\hbox{\rm d}}}

\def\C{{\hbox{\bb C}}}

\def\hop{\vskip 4pt plus 2pt} 
\def\hophop{\vskip 7pt}

\begin{document} 

% Double spacing

\begin{titlepage} 
\begin{center} 
\mbox{ }\vspace{2in} \\ 
\huge\bfseries 
Completions of normed algebras of differentiable functions \vspace{1in} \\ 
\normalsize\normalfont 
William J. Bland and Joel F. Feinstein$^*$\vspace{1in} \\ 
\large 
This paper includes work from 
the doctoral\\ thesis of the first author \cite{Bland_Thesis} 
\end{center} 
\end{titlepage} 

%\tableofcontents
%\chapter{}

\abstract{In this paper we look at normed spaces of differentiable functions 
on compact plane sets, including the spaces of infinitely differentiable 
functions considered by Dales and Davie in \cite{Dales+Davie}. For many 
compact plane sets the classical definitions give rise to incomplete spaces. 
We introduce an alternative definition of differentiability which allows us to 
describe the completions of these spaces. We also consider some associated 
problems of polynomial and rational approximation.} 

\section{Introduction} 

In this paper we shall investigate problems concerning the normed spaces of 
differentiable functions on compact plane sets which were originally studied 
by 
Dales, Davie and McClure with particular reference to the 
Banach algebra case in \cite{Dales+Davie} and \cite{Dales+McClure}. These 
spaces are the $D(X,M)$ spaces where $X$ is a perfect 
compact plane set and $M$ is a 
sequence of positive real numbers. We shall give the definitions of these 
and some of the related spaces that we wish to study in Section 2.

There are many interesting problems concerning these spaces. 
One strand of their study concerns problems of approximation of functions 
by means of polynomials or rational functions. 
For example, in \cite{Dales+McClure} Dales and McClure proved that, 
if $X$ is the closed unit disk, then the polynomials are always dense in 
$D(X,M)$. 
In the one-dimensional case, some results on polynomial approximation were 
obtained by O'Farrell in \cite{OFar} (see also Theorem 4.4.15 
of \cite{Garthsbook}). It was shown under very mild conditions 
on $M$ that the 
polynomials are dense in $D(I,M)$, where $I$ is a closed interval. 
Further results on holomorphic and polynomial approximation, and related problems 
concerning extensions of functions in these spaces, were obtained in 
\cite{FLO}. There are still many fascinating open problems in this area. 

Another set of problems concerns the completeness of these spaces. 
Dales and Davie (\cite{Dales+Davie}) gave some 
conditions on $X$ which guaranteed the completeness of the $D(X,M)$ spaces. 
It was also noted in \cite{Ali_Thesis} that a union of any finite number of 
sets for which these spaces are complete gives another such set. Some further results 
on this are given in \cite{Honary+Mahyar1}.
We shall 
discuss this problem further in Section 2. We 
prove that for any perfect, compact plane set which has infinitely many 
components all of these spaces are 
incomplete.
We also give an 
example of an $X$ which is a rectifiable Jordan arc and yet the 
normed algebra of once continuously differentiable functions on $X$ (with the 
classical definition) is incomplete. 

In the setting of normed or Banach algebras, work on these spaces includes the 
study of their endomorphisms. 
In \cite{Kamowitz}, Kamowitz came close to classifying all the endomorphisms 
of $D(I,M)$, 
in terms of self-maps of the interval $I$. The second author and Kamowitz made 
further progress on the remaining problems, and investigated more general 
compact plane sets in \cite{Feinstein+Kamowitz}.
We will not investigate these problems in this paper, but it is worth noting 
that some of the results of \cite{Feinstein+Kamowitz} depend on the completeness 
of the $D(X,M)$ spaces. (See also \cite{Behrouzi} for some work on homomorphisms
between these spaces.)

In Section 3 we discuss various matters related to rectifiable arcs, including 
the Fundamental Theorem of Calculus for rectifiable paths, and the conditions 
of uniform regularity and pointwise regularity for compact plane sets. 
These 
latter conditions are sufficient to imply completeness of all the normed 
spaces defined in Section 2 (\cite{Dales+Davie}, \cite{Honary+Mahyar1}). 
In the case where the spaces are incomplete, it becomes important to 
investigate their completions. 
We do this in Section 4. 
There we determine the completions of the normed spaces above, at least 
for compact plane sets $X$ such that the union of the images of all the 
injective, rectifiable arcs in $X$ is a dense subset of $X$. 
In this setting we define a less restrictive 
notion of differentiation which ensures that the spaces we end up with are 
complete. The original versions of the spaces embed isometrically in our new 
versions, so 
the completions of the original spaces are simply their closure in the new 
spaces. Where the algebras considered in \cite{Feinstein+Kamowitz} were 
incomplete, the new versions are complete and all the arguments of Feinstein 
and Kamowitz remain valid in the new setting. This suggests that the new 
algebras may, in fact, be the correct place to study endomorphisms. 

In Section 5 we investigate two related problems. 
For which compact plane sets are the new spaces constructed in Section 4 
the same as the original spaces as defined in Section 2? 
For which compact plane sets are the 
original spaces dense in the new spaces? We also obtain some related 
polynomial and rational approximation results for these spaces. 
For some work on identifying the maximal ideal 
spaces in the original setting
and on polynomial and rational/holomorphic approximation in some related 
spaces 
see, for example,
\cite{FLO},
\cite{Honary}, \cite{Honary+Mahyar1},
\cite{Honary+Mahyar2},
\cite{Jarosz}, 
\cite{OFar} and \cite{Verdera}.

We conclude, in Section 6, with some open problems. 

\section{Introductory concepts and results} 

We begin with some standard notation, definitions and results. 
Let $X$ be a compact plane set. 
We denote the set of all continuous, complex-valued functions on
$X$ by $C(X)$. For $f \in C(X)$ we denote the uniform norm of $f$ by
$\|f\|_\infty$. More generally we denote the uniform norm of $f$ on a closed subset $E$ of
$X$ by $\|f\|_E$.

\begin{definition} 
Let $X$ be a perfect, compact plane set $X$. 
We say that a complex-valued function $f$ defined on $X$ 
is {\it complex-differentiable} 
at a point $a\in X$ if 
the limit 
$$ f^\prime(a) = \lim_{z\to a, \ z\in X} 
{f(z)-f(a)\over z-a} $$ 
exists. We call $f^\prime(a)$ the {\it complex derivative} of $f$ at $a$. 
Using this concept of derivative, we define the terms 
{\it complex--differentiable on} $X$, 
{\it continuously complex--differentiable on} $X$, 
and {\it infinitely complex--differentiable on} $X$ 
in the obvious way. We denote the $n$-th complex derivative 
of $f$ at $a$ by $f^{(n)}(a)$, and we denote the set of infinitely 
complex-differentiable functions on $X$ by $D^\infty(X)$. 
We denote the set of continuously complex-differentiable functions 
on $X$ by $\stddiff{X}$. More generally, we define the corresponding algebras 
of $n$-times continuously differentiable functions, $D^n(X)$, again in the 
obvious way. 

Let $(M_n)$ be a sequence of positive real numbers. 
We define the space $\dis D(X,M)=\{f \in D^{\infty}(X): 
\|f\|=\sum_{n=0}^{\infty} 
\frac{\|f^{(n)}\|_{\infty}}{M_n} < \infty\}.$ With pointwise addition, 
$D(X,M)$ is a normed space which 
is not necessarily complete. 

If further the sequence $M_n$ satisfies 
$M_0=1$ and, for all non-negative integers $m$, $n$, we have 
$$\dis \left(\begin{array}{c}m+n\\n \end{array} \right) \leq 
\frac{M_{m+n}}{M_mM_n}$$ 
then $D(X,M)$ is a normed algebra with pointwise multiplication. 

In \cite{Dales+Davie}, Dales and Davie used this class of algebras to give an 
example of a commutative semisimple 
Banach algebra for which the peak points are of first category 
in the Silov boundary, and an example of a commutative 
semisimple Banach algebra $B$ and a discontinuous function $F$ acting on $B.$ 

Clearly the restrictions of all (analytic) polynomials to $X$ 
belong to all of the $D(X,M)$ spaces. It was further proved 
in \cite{Dales+Davie} that 
the algebra $D(X,M)$ includes 
all of the rational functions with poles off $X$ if and only if 
$$\lim_{n\to\infty} {\left({n!}\over{M_n}\right)}^{1\over n} = 0. \eqno(1)$$ 
We say that $(M_n)$ is a {\it nonanalytic sequence} if (1) holds 
\cite{Dales+Davie}. 
\end{definition}

Each of the spaces $D^n(X)$ is a normed algebra, using the norm 

\[\|f\|=\sum_{k=0}^n{\frac{\|f^{(k)}\|_\infty}{k!}}.\] 

These spaces are often incomplete, even for fairly nice $X$: we give
some examples of this below. However, for a given $X$, the completeness of 
$\stddiff{X}$ implies the completeness of all of the others. This follows from
the following result. 

\begin{theorem} 
Let $X$ be a perfect, compact plane set and let $r$ be a positive integer.
Suppose that $D^r(X)$ is 
complete. Then, for all integers $n\geq r$, $D^n(X)$ is complete and, 
for every sequence $M$ of positive real numbers, $D(X,M)$ is 
complete. 
\end{theorem} 

\jremth{induce_completeness}
\Proof 
We give the proof for $D(X,M)$. The proof for $D^n(X)$ is similar but slightly 
easier. 
Let $f_m$ be a Cauchy sequence in $D(X,M)$. It is clear that, for all 
non-negative integers $k$ the sequence $f_m^{(k)}$ is Cauchy in 
$D^r(X)$ and so converges in $D^r(X)$ to a function $g_m$, say. By 
definition of the norm on $D^r(X)$, we see that the sequence 
$(f_m^{(k)})'$ converges uniformly to $g_m '$ on $X$. However, we also know 
that $(f_m^{(k)})'=f^{(m+1)}$ converges to $g_{m+1}$ and 
so we have $g_m ' = g_{m+1}$. 
The remainder of the proof is a standard functional analysis argument 
showing that $g_0 \in D(X,M)$ and that the sequence $f_m$ converges 
in $D(X,M)$ to $g_0$: we omit the details. 

\QED 

We now prove that if $X$ has infinitely many 
components then $D(X,M)$ is incomplete, and hence all of
the spaces $D^n(X)$ are 
incomplete. In the proof, and throughout the rest of this paper, 
we will frequently refer to sets which are both open and closed, and it will be
convenient to call such sets {\it clopen} sets.

\begin{theorem} 
Let $X \subseteq \Complex$ be a compact, perfect set, which has 
infinitely many components and let $M$ be any
sequence of positive real numbers. 
Then all of the spaces $D^n(X)$ and $D(X,M)$ are incomplete 
\end{theorem} 
\jremth{Many_Components} 

\Proof 
By Theorem \jrecth{induce_completeness}
it is sufficient to prove the result for $D(X,M)$.
(The proof given below is, anyway, valid in all cases.)
We are given that $X$ has infinitely many connected components. 
Set $E_{0} = X$. Then $E_{0}$ can be written as $E_{0} = E_{1} \cup F_{1}$ 
where $E_{1}$ and $F_{1}$ are nonempty, disjoint, clopen subsets of $E_{0}$ and 
$E_{1}$ has infinitely many components. 

Similarly we can write $E_{1} = E_{2} \cup F_{2}$ where 
$E_{2}$ and $F_{2}$ are nonempty disjoint clopen subsets of $E_{1}$ and 
$E_{2}$ has infinitely many components. 

Clearly we can continue in this way to form sequences $(E_{n})$ and $(F_{n})$. 
For each $n \in \Natural$, choose a point $z_{n} \in F_{n}$. Then the sequence $z_n$ 
has a convergent subsequence, $z_{n_{k}}$. 
Say $z_{n_{k}} \rightarrow z_{0}$ as $n \rightarrow \infty$. 
Now, we cannot have $z_{0} \in F_{n_{k}}$ for any $k \in \Natural$ since the 
sets 
$F_{n_{k}}$ are open and pairwise disjoint. 

Define $f \in C(X)$ by 
\[ 
f(z) = \left\{ 
\begin{array}{ll} 
z_{n_{k}} & \mbox{ for } z \in F_{n_{k}} \\ 
z_{0} & \mbox{ for } z \in X \setminus \bigcup_{k=1}^{\infty} F_{n_{k}} 
\end{array} 
\right. 
\] 
Then $f$ is constant on each of the clopen sets $F_{n_{k}}$ and so has
derivative $0$ on their union. Thus if $f$ was in $\stddiff{X}$
we would also have $f'(z_0)=0$.
However, for all $k$, we have
$\frac{f(z_{n_{k}})-f(z_{0})}{z_{n_{k}}-z_{0}} 
 = 1$,
and so $f$ is not in 
$\stddiff{X}$. 
Finally, note that
there is an obvious sequence $(f_{i}) \subseteq D(X,M)$ such that
$f_{i} \rightarrow f$ uniformly on $X$:
for $i \in \Natural$, define $f_{i} \in \stddiff{X}$ by 
\[ 
f_{i}(z) = \left\{ 
\begin{array}{ll} 
z_{n_{k}} & \mbox{ if } z \in F_{n_{k}} \mbox{ and } k \leq i \\ 
z_{0} & \mbox{ for } z \in X \setminus \bigcup_{k=1}^{i} F_{n_{k}} 
\end{array} 
\right. 
\] 

It is easy to see that $f_i ' = 0$ for all $i$ and that $(f_{i})$
is a Cauchy sequence in 
$D(X,M)$.
Since $f$ is not even in $\stddiff{X}$,
$D(X,M)$ is incomplete. 

\QED 

\newcounter{pwreg_complete_parcon} 
\remembertheorem{pwreg_complete_parcon} 

The completeness of
$\stddiff{X}$ is far from being a topological property of $X$:
we conclude this section with an example
where $X$ is the image of a rectifiable Jordan arc in the plane
and yet $\stddiff{X}$ is 
incomplete. (We will look at rectifiable curves in more detail later 
in this paper.) 
\begin{example} 
Set $z_{n} = 2^{-2n}$ and $w_{n} = 2^{-2n}+2^{-n}i$. 
We glue together the origin and the following paths ($\gamma_{n}$ for $n \in 
\Natural$): 

{ 
\psfrag{xnp1}{$2^{-2(n+1)}$} 
\psfrag{xn}{$2^{-2n}$} 
\psfrag{yn}{$2^{-n}$} 
\psfrag{zn}{$z_{n}$} 
\psfrag{wn}{$w_{n}$} 
\psfrag{origin}{Origin} 
\includegraphics{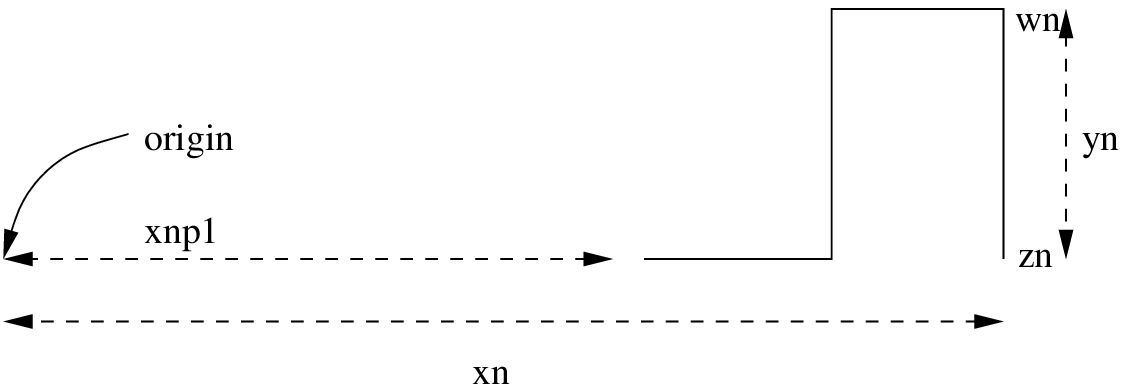} 
} 

The resulting path $\gamma$ can be parametrised by its arc-length. 
It is clear that $\gamma$ is a rectifiable Jordan arc.

The exact position on the x-axis of the leftmost vertical line forming 
$\gamma_{n}$ is 
irrelevant to the working of this example, so long as it lies (strictly) 
between $2^{-2(n+1)}$ and $2^{-2n}$. 
\end{example} 

\newcounter{example_path_incpl} 
\remembertheorem{example_path_incpl} 
\jremth{example_path_incpl} 

\begin{theorem} 
Let $X$ be the image of the path $\gamma$ in the previous example. 
Then $\stddiff{X}$ is incomplete. 
\end{theorem} 
\Proof 
Define $f \in C(X)$ by 
\[ 
f(0) = 0 
\] 
and 
\[ 
f(x+yi) = \left\{ 
\begin{array}{ll} 
-3 \cdot 2^{n-1}y^{3} + \frac{9y^{2}}{4} + 2^{-2n} & \mbox{if $x = 2^{-2n}$} 
\\ 
-3 \cdot 2^{n-1} \cdot 2^{-n} + \frac{9 \cdot 2^{-2n}}{4} + 2^{-2n} & 
\mbox{otherwise} 
\end{array} 
\right. 
\] 
where $x+yi \in \Image{\gamma_{n}}$. 

It is straightforward to check that the following conditions hold:
\begin{enumerate} 
\item $f$ is constant everywhere on $\Image{\gamma_{n}}$ except on the line 
joining the points $z_{n}$ and $w_{n}$; 
\item $f$ is continuous on the whole of $\Image{\gamma}$, and is continuously 
differentiable on each $\Image{\gamma_{n}}$;
\item $f'(z_{n}) = f'(w_{n}) = 0$; 
\item $f(z_{n}) = z_{n}$; 
\item $f(w_{n}) = z_{n+1}$; 
\item as $n \rightarrow \infty$, the supremum of $|f'|$ on the image of 
$\gamma_{n}$ converges to zero. 
\end{enumerate} 
Thus we have 
\[ 
\lim_{n \rightarrow \infty} \frac{f(z_{n})-f(0)}{z_{n}-0} 
= 
\lim_{n \rightarrow \infty} \frac{z_{n}-0}{z_{n}-0} 
= 
1 
\neq 
0 
= 
\lim_{n \rightarrow \infty} f'(z_{n}) 
\] 
and so $f \notin \stddiff{X}$. 

However, there is an obvious Cauchy sequence $(f_{n})$ of functions in 
$\stddiff{X}$ 
with $\|f-f_{n}\|_{\infty} \rightarrow 0$. We simply define $f_{n}(z)$ to be 
equal 
to $f(z)$ when the real part of $z$ is larger than $z_{n}$, and constantly 
equal 
to $f(z_{n})$ otherwise. 
Thus $\stddiff{X}$ is incomplete. 

\QED

\section{Rectifiable paths and regularity conditions for compact plane sets}

In this section we discuss families of rectifiable curves and 
some related conditions. 
We will assume that the reader is familiar with the elementary results and 
definitions concerning rectifiable paths including integration of continuous, 
complex-valued functions along rectifiable curves. For more details see, for 
example, Chapter 6 of \cite{Apostol}. 
\hophop 
\begin{definition} 

The paths considered in this paper will be 
continuous, complex-valued functions defined on {\it non-degenerate} closed 
intervals $[a,b]$: 
$a$, $b$ will be real numbers with $a<b$. 
When $\gamma$ is such a path, we say that it is a path {\it from 
$\gamma(a)$ to $\gamma(b)$} and we denote these {\it endpoints} of the path by 
$\gamma^-$ and $\gamma^ +$. 
\hop 
Given $X\subseteq \C$, a {\it path in $X$} is a path whose image is 
a subset of $X$. (A Jordan arc in $X$ is, of course, 
simply an injective path in $X$.)

%We denote the union of the images of all paths in $X$ by %$\allpaths{X}$.

The length of a rectifiable path $\gamma$ will be denoted by $|\gamma|$.

\end{definition} 

We recall the following elementary connection between piecewise smooth paths, 
rectifiability and integration. 

\begin{proposition}[see \cite{Conway}, pp.58-62] 
Let $X$ be a compact subset of $\Complex$ and $\gamma : [a,b] \rightarrow X$ 
be a 
piecewise smooth path in $X$. Then: 
\begin{enumerate} 
\item $\gamma$ is rectifiable; 
\item $\Length{\gamma} = \int_{a}^{b} |\gamma'(t)| \wrt{t}$; 
\item $\int_{\gamma} f(z) \wrt{z} = \int_{a}^{b} f(\gamma(t)) \gamma'(t) 
\wrt{t}$ 
\mbox{ for any $f \in C(X)$}. 
\end{enumerate} 
\end{proposition} 

The next result is an analogue of the 
Fundamental Theorem of Calculus. We have not been able to find a proof of 
it in the literature, although a similar theorem is given in \cite{Conway} 
(Theorem 1.18, p.65). 
However, the functions in that theorem are defined on open subsets of 
$\Complex$, whereas we need the same result for functions defined only on 
images of rectifiable paths. 
Elegant proofs of this general result using 
the method of repeated bisection have been given to us by G. R. Allan, T.W. 
K\"{o}rner 
and W.K. Hayman.
The proof provided by Allan may be found in full in
\cite{Bland_Thesis}.
\begin{theorem} {\bf (Fundamental theorem of calculus for rectifiable paths)} 
Let $\gamma : [a,b] \rightarrow \Complex$ be a rectifiable path with endpoints 
$z_{1}$ and $z_{2}$. Then for every $f \in \stddiff{\Image{\gamma}}$ we have 
\[\int_{\gamma} f'(z) \wrt{z} = f(z_{2}) - f(z_{1}).\] 
\end{theorem}

We now discuss, in terms of rectifiable paths, some standard conditions a compact plane 
set $X$ may satisfy which are sufficient to ensure the completeness of 
$\stddiff{X}$ (and hence of all the other spaces 
defined in Section 2). 
In \cite{Ali_Thesis} it was shown that the collection of sets
$X$ for which  $\stddiff{X}$ is complete is closed under
finite union.
(In fact the result is only stated there for the 
$D(X,M)$ spaces, but the 
proof for the other spaces is the same). 

\begin{definition} 
Let $X \subseteq \Complex$ be compact. We say $X$ is \emph{regular} at a 
point $z \in X$ if there is a constant $k_{z} > 0$ such that, 
for every $w \in X$ there is a path $\gamma : [a,b] \rightarrow X$ 
with $\gamma(a) = z$, $\gamma(b) = w$ and $\Length{\gamma} \leq k_{z}|z-w|$. 

We say $X$ is \emph{pointwise regular} if $X$ has more than one point and 
$X$ is regular at every point $z \in X$. (In \cite{Honary+Mahyar1} such a set is simply said to
be {\it regular}.)
If, further, there is one constant $k > 0$ such that, for all $z$ and $w$ in 
$X$, 
there is a path $\gamma : [a,b] \rightarrow X$ with $\gamma(a) = z$, 
$\gamma(b) = w$ and $\Length{\gamma} \leq k|z-w|$ then we say that $X$ is 
\emph{uniformly regular}. 
\end{definition} 

Clearly all pointwise and uniformly regular sets are perfect and 
path-connected. 
For points $z$ and $w$ in a set $X \subseteq \Complex$, we will define 
\[ 
d(z,w) = 
\inf \{|\gamma| : \gamma \mbox{ is a rectifiable path from $z$ to $w$ in 
$X$}\}. 
\] 

Dales and Davie showed that $\stddiff{X}$ is complete whenever
$X$ is a finite union of 
uniformly regular sets.
However, as observed in \cite{Honary+Mahyar1}, the proof given 
in \cite{Dales+Davie} 
is equally valid for pointwise regular sets. Thus
$\stddiff{X}$ is complete whenever
$X$ is a finite union of pointwise regular sets.
We will give another proof of this fact in the next section 
where we investigate the completions of these normed spaces. 

We now note that for a compact plane set $X$ to satisfy one of these 
two regularity conditions it is sufficient (though of course not necessary) 
for the boundary to satisfy the same condition. 

\begin{theorem} 
Let $X \subseteq \Complex$ be compact. If $\Boundary{X}$ is uniformly 
regular then $X$ is uniformly regular. 
\end{theorem} 
\Proof 
We know that there is a constant $k > 0$ such that $d(z,w) \leq k|z-w|$ 
for all $z,w \in \Boundary{X}$. 
Choose $z_{1}, z_{2} \in X$. 
If the line-segment connecting $z_{1}$ and $z_{2}$ is contained in $X$ 
then $d(z_{1}, z_{2}) = |z_{1}-z_{2}|$ and we are done. 
Otherwise, the line-segment must 
intersect 
$\Boundary{X}$ at at least two points. 
In this case,
let $Z$ be the set of points of intersection, and (for $i=1,2$)
let $w_{i}$ be the closest point of $Z$ to $z_{i}$.

We know that there is a path $\gamma$ in $\Boundary{X}$ between 
$w_{1}$ and $w_{2}$ such that $|\gamma| \leq k|w_{1}-w_{2}|$. 
We have 
\begin{eqnarray*} 
d(z_{1},z_{2}) & \leq & |\gamma| + |z_{2}-w_{2}| + |w_{1}-z_{1}| \\ 
& \leq & k|w_{1}-w_{2}| + |z_{1}-w_{1}| + |z_{2}-w_{2}| \leq  (k+1)|z_{1}-z_{2}|. 
\end{eqnarray*} 

\QED 

The proof of the same theorem for pointwise regularity requires a little more 
thought, but is essentially the same. 

\begin{theorem} 
Let $X \subseteq \Complex$ be compact. If $\Boundary{X}$ is pointwise 
regular then $X$ is pointwise regular. 
\end{theorem} 
\Proof 
Choose $z_{1} \in X$. Choose a point $w_{1} \in \Boundary{X}$ such that 
no other point in $\Boundary{X}$ is closer to $z_{1}$ than $w_{1}$. 
We know that there is a constant $k_{w_{1}} > 0$ such that for any 
$w_{2} \in \Boundary{X}$ we have $d(w_{1}, w_{2}) \leq 
k_{w_{1}}|w_{1}-w_{2}|$. 
Set $c_{z_{1}} = 2+3k_{w_{1}}$. 

Choose $z_{2} \in X$. Again, choose a point $w_{2} \in \Boundary{X}$ such that 
no other point in $\Boundary{X}$ is closer to $z_{2}$ than $w_{2}$. 

If the line-segment connecting $z_{1}$ to $z_{2}$ is contained in $X$ then we 
are done. Assume otherwise. 
Then we have 
\[|z_{1}-w_{1}| \leq |z_{1}-z_{2}|\] 
and 
\[|z_{2}-w_{2}| \leq |z_{1}-z_{2}|\] 
Thus 
\begin{eqnarray*} 
d(z_{1}, z_{2}) & \leq & k_{w_{1}}|w_{1}-w_{2}| + |z_{1}-w_{1}| + 
|z_{2}-w_{2}| \\ 
& \leq & k_{w_{1}}|w_{1}-w_{2}| + 2|z_{1}-z_{2}| 
\end{eqnarray*} 
Now 
\[ 
|w_{1}-w_{2}| = |(w_{1}-z_{1})+(z_{1}-z_{2})+(z_{2}-w_{2})| \leq 
3|z_{1}-z_{2}| 
\] 
and so 
\[ 
d(z_{1}, z_{2}) \leq c_{z_{1}}|z_{1}-z_{2}| 
\] 
as required. 

\QED 

The following elementary lemma will be useful later. 

\begin{lemma} 
Each component of a finite union of pointwise regular sets is pointwise 
regular. 
\end{lemma} 
\Proof 
Let $X$ and $Y$ be pointwise regular compact plane sets such that $X \cap Y 
\neq \emptyset$. It is clear that we just need to show that $X \cup Y$ is 
pointwise regular and then the proof follows. 

Choose $z \in X \cup Y$. If $z \in X \cap Y$ then the result is clear. 
Assume without loss of generality that $z \in X \setminus Y$. 
Then $\mbox{dist}(z,Y) > 0$, where $\mbox{dist}(z,Y)$ is the usual Euclidean 
point-set distance between $z$ and $Y$. Using this and the pointwise regularity of $X$
it is now elementary 
to show that there is a $c>0$ such that for all $w \in X \cup Y$ we have 
$d(z,w) \leq c|z-w|$. 
\QED 

\newcounter{fin_pr} 
\remembertheorem{fin_pr} 

The property of pointwise regularity is a local property in the following 
sense: say a set $X$ is {\it locally pointwise regular} if each point in $X$ 
has a pointwise regular compact neighbourhood in $X$. 
The following result is now an immediate consequence of compactness. 
\begin{theorem} 
Let $X \subseteq \Complex$ be a locally pointwise regular, compact 
plane set. 
Then $X$ is a finite union of pointwise regular sets. 
\end{theorem} 

We finish this section by noting that there are examples of rectifiable paths 
in the 
complex plane whose images are not pointwise regular. 
For example, the fact that the path $\gamma$ 
in Example \jrecth{example_path_incpl} has $\stddiff{\Image{\gamma}}$ 
incomplete implies that $\gamma$ cannot be pointwise regular. 

We are now ready to introduce the new normed spaces that we wish to study. 

\section{The $\paths$-differentiation spaces} 

In this section we investigate the completions of the normed spaces 
considered above by weakening the differentiability requirement on the 
functions. 

One well-known, related class of Banach spaces is to look at analytic 
functions on an open subset $U$ of $\Complex$ 
with some specified number of the function's derivatives 
being bounded. This gives a set of Banach spaces corresponding to the spaces 
$D^n(X)$ above. Indeed, when $X$ is the closure of $U$, the spaces $D^n(X)$ 
embed isometrically in these new, complete, spaces. 
A similar construction provides complete versions of the spaces $D(X,M)$. 
However, these constructions are only helpful for compact spaces $X$ where the 
interior of $X$ is dense in $X$. This is too restrictive for our purposes. 
Instead, we will mostly work with the larger class of compact plane sets $X$ 
for which the union of the images of all rectifiable Jordan arcs
in $X$ is dense in $X$. We will then use appropriate sets of Jordan
arcs to define our notion of derivative.

\begin{definition} 
Let $X \subseteq \Complex$ be compact and $\paths$ be a set of paths in $X$. 
We say $\paths$ is \emph{useful} if the following conditions are satisfied.
\begin{enumerate} 
\item Every path in $\paths$ is a rectifiable Jordan arc.
\item If $\gamma \in \paths$ is defined on $[a,b]$, then the restriction 
of $\gamma$ to $[c,d]$ is in $\paths$ whenever $[c,d] \subseteq [a,b]$ and $c 
< d$. 
\end{enumerate} 
We write 
\[\paths(X) = \bigcup_{\gamma \in \paths} \Image{\gamma}\] 
and for $z, w \in X$ we set 
\[\paths(z,w) = \{ \gamma \in \paths : \gamma^{-} = z, \gamma^{+} = w \}\] 
(recall that $\gamma^{-}$ and $\gamma^{+}$ are the endpoints of $\gamma$). 
\end{definition} 

Clearly the sets of rectifiable Jordan arcs and smooth, rectifiable Jordan
arcs in $X$ 
are both useful. Also, for any $L>0$, the set of rectifiable Jordan
arcs 
in $X$ with length $\leq L$ is useful.

We are now ready to define the notion of differentiability associated
with with a set of rectifiable paths.

\begin{definition} 
Let $X \subseteq \Complex$ be compact and $\paths$ be a set of rectifiable 
paths 
in $X$. For $f \in C(X)$, we say $g \in C(X)$ is an \emph{$\paths$-derivative} 
of $f$ 
if, for all $z_{1},z_{2} \in X$ and every $\gamma \in \paths(z_{1},z_{2})$ 
we have \[\int_{\gamma} g(z) \wrt{z} = f(z_{2})-f(z_{1}).\] 
\end{definition} 

% In future, could require g to be continuous 
% along each path gamma rather than being
% in C(X).

We will mostly restrict attention to the case where 
$\paths$ is a useful set of paths. 

\begin{definition} 
Let $X \subseteq \Complex$ be compact and $\paths$ be a set of rectifiable 
paths in $X$. Define 
\[ \diff{\paths}{X} = \{ f \in C(X) : f \mbox{ has an $\paths$-derivative in 
$C(X)$} \}.\] 
\end{definition} 

Clearly we would not expect $\paths$-derivatives to be unique. We 
will see below, however, that their restriction to 
$\Closure{\paths(X)}$ is unique. 

The following theorem is the $\paths$-derivative analogue of a standard result 
of elementary real analysis. 

\begin{theorem} 
Let $X$ be a compact plane set and let $\paths$ be a useful set of paths in 
$X$. Let $f_n$, $g_n$ be uniformly convergent sequences in $C(X)$ with limits 
$f$, $g$ respectively. Suppose that, for all 
$n$, $g_n$ is an $\paths$-derivative of $f_n$. Then $g$ is an 
$\paths$-derivative of $f$. 
\end{theorem} 
\newcounter{lim_deriv} 
\remembertheorem{lim_deriv} 

\Proof 
This is essentially immediate from the definitions.
\QED 

As we have already seen, the analogous statement for the original notion of 
differentiation is false: this is the reason why the 
spaces $\stddiff{X}$ are often incomplete. 

Before going any further, we 
deal with the issue of ``piecewise'' curves. 
For every set $\paths$ of paths in $X$, there is a corresponding set 
$\paths_{PW}$ 
of paths that are ``piecewise-$\paths$'' paths. In other words, each path in 
$\paths_{PW}$ 
consists of finitely many paths in $\paths$ that are joined together at their 
endpoints. 
The question is, do $\paths$ and $\paths_{PW}$ lead to different theories 
of differentiation? 

\begin{theorem} 
Let $X \subseteq \Complex$ be compact and $\paths$ be a useful set of paths 
in $X$. Let $\paths_{PW}$ be the piecewise version of $\paths$, as described 
above. Then $\diff{\paths}{X} = \diff{\paths_{PW}}{X}$. 
\end{theorem} 
\Proof 
Again this is an elementary consequence of the definitions.
\QED 

In view of this result we can now take $\paths$ to be (for example) 
either 
the set of smooth Jordan arcs in $X$, or the set of piecewise 
smooth Jordan arcs in $X$; 
it does not make any difference to the resulting object $\diff{\paths}{X}$. 
Also we may assume that the lengths of the paths in $\paths$
are bounded: for example, for each $L>0$, the same theory is obtained by using the set of 
all rectifiable Jordan arcs in $X$ as is obtained 
by using the set of those  
rectifiable Jordan arcs in $X$ whose length is at most $L$.
(Every rectifiable curve is \lq piecewise short'.)

Note that we always have 
$\Closure{\paths(X)} = \Closure{\paths_{PW}(X)}$. 
Also, the set $\paths_{PW}$ is useful if $\paths$ is useful 
(the converse is, however, not true). 

We will prove that $\diff{\paths}{X}$ is always a Banach algebra. As part of this we need to
check that 
$\paths$-derivatives behave in the way we expect with regard to sums,
scalar multiples and products.

\begin{theorem} 
Let $X \subseteq \Complex$ be compact and $\paths$ be a useful set of paths 
in $X$. Let $f_{1}, f_{2} \in C(X)$ and $\lambda \in \Complex$. 
If $g_{1}, g_{2} \in C(X)$ are $\paths$-derivatives of $f_{1}$ and $f_{2}$ 
respectively 
then $g_{1}+\lambda g_{2}$ is an $\paths$-derivative of $f_{1}+\lambda f_{2}$. 
\end{theorem} 
\Proof 
Set $f = f_{1} + \lambda f_{2}$ and $g = g_{1} + \lambda g_{2}$. Clearly $g 
\in C(X)$. 
Now choose $z_{1}, z_{2} \in X$ and $\gamma \in \paths(z_{1},z_{2})$. 
We have: 
\begin{eqnarray*} 
\int_{\gamma} g(z) \wrt{z} 
& = & \int_{\gamma} (g_{1}(z) + \lambda g_{2}(z)) \wrt{z} \\ 
& = & \int_{\gamma} g_{1}(z) \wrt{z} + \int_{\gamma} \lambda g_{2}(z) \wrt{z} 
\\ 
& = & f_{1}(z_{2}) - f_{1}(z_{1}) + \lambda (f_{2}(z_{2}) - f_{2}(z_{1})) \\ 
& = & f(z_{2}) - f(z_{1}) 
\end{eqnarray*} 
Thus $g$ is an $\paths$-derivative of $f$. 

\QED 

\newcounter{difflinear} 
\remembertheorem{difflinear} 

\begin{corollary} 
Let $X \subseteq \Complex$ be compact and $\paths$ be a useful set of paths 
in $X$. Then $\diff{\paths}{X}$ is a vector space. 
\end{corollary} 

We now look at multiplication of functions in $\diff{\paths}{X}$. 
First we note an elementary lemma concerning polynomials. 

\begin{lemma} 
Let $X \subseteq \Complex$ be compact and $\gamma$ be an injective rectifiable 
path 
whose image is contained in $X$. Set $\paths = \{\gamma\}$. 
Then for any polynomials $p_{1}$ and $p_{2}$ defined on $X$, the function 
$p_{1}'p_{2} + p_{1}p_{2}'$ is an $\paths$-derivative of $p_{1}p_{2}$. 
\end{lemma} 
\Proof 
We know that $p_{1}$, $p_{2}$ and $p_{1}p_{2}$ are all complex-differentiable 
on any complex plane set, and $(p_{1}p_{2})' = p_{1}'p_{2} + p_{1}p_{2}'$. 
The result now follows from the Fundamental Theorem of Calculus for 
rectifiable 
paths (or indeed the special case of this theorem for polynomial functions). 

\QED 

\begin{theorem} 
Let $X \subseteq \Complex$ be compact and $\gamma$ be an injective rectifiable 
path whose image is contained in $X$. Let $\paths$ be the set of all 
subpaths of $\gamma$ (including $\gamma$ itself). 
Then for any functions $f_{1}, f_{2} \in \diff{\paths}{X}$ with 
$\paths$-derivatives $g_{1}$ and $g_{2}$ respectively, the function 
$g_{1}f_{2} + f_{1}g_{2}$ is an $\paths$-derivative of $f_{1}f_{2}$. 
\end{theorem} 
\Proof 
Set $Y = \Image{\gamma}$. 

Note that, since $\gamma$ is injective, 
$\C \setminus \Image{\gamma}$ must be connected.
Hence by Mergelyan's (or Lavrentiev's) theorem we can choose two sequences of 
analytic polynomials $p_n$, $q_n$ converging uniformly on $Y$ to $g_1$, $g_2$ 
respectively. Now (antidifferentiating) choose analytic polynomials 
$P_n$, $Q_n$ such that $P_n ' = p_n$, $Q_n '=q_n$, 
$P_n(\gamma^-)=f_1(\gamma^-)$ and $Q_n(\gamma^-)=f_2(\gamma^-)$. 

The Fundamental Theorem of Calculus for rectifiable paths tells us that $p_n$ 
is an $\paths$-derivative of $P_n$ and similarly for 
$q_n$ and $Q_n$. It now follows easily that $P_n$, $Q_n$ converge 
uniformly on $Y$ to $f_1$, $f_2$ respectively. 

By the preceding lemma we know that 
$P_n q_n + p_n Q_n$ is an $\paths$-derivative of $P_n Q_n$. 
Taking uniform limits, and applying \recalltheorem{lim_deriv}, 
we see that $f_1 g_2 + f_2 g_1$ is an $\paths$-derivative 
of $f_1 f_2$, as required. 

\QED 

\begin{corollary} 
Let $X \subseteq \Complex$ be compact and $\paths$ be a useful set of paths 
in $X$. 
Then $\diff{\paths}{X}$ is an algebra 
\end{corollary} 
\Proof 
Clear. 
\QED

We now wish to establish the extent to which 
$\paths$-derivatives are unique. 
We start with a simple lemma. 

\begin{lemma} 
Let $X \subseteq \Complex$ be compact and $\gamma$ be a rectifiable path in 
$X$ defined 
on $[a,b]$ and with $\gamma(a) \neq \gamma(b)$. 
Then there exists a constant $k > 0$ and a sequence $(\gamma_{n})$ of 
sub-paths of $\gamma$ defined on subintervals $[a_n,b_n]$ of $[a,b]$, with 
$|\gamma_n|=2^{-(n-1)}|\gamma|$ and such that 
for each $n \in \Natural$, $|\gamma_{n}| < k|\gamma_{n}(b_{n}) - 
\gamma_{n}(a_{n})|$. 
\end{lemma} 
\Proof 
Clearly $|\gamma| < k|\gamma(b)-\gamma(a)|$ for some $k > 0$. This will be our 
$k$. We set $a_1=a$, $b_1=b$ and $\gamma_1=\gamma$. 
Now suppose that $n>1$ and that we have constructed the sequence of sub-paths 
up to and including the path 
$\gamma_{n-1}$, which is defined on $[a_{n-1},b_{n-1}]$. 
Choose $c \in (a_{n-1},b_{n-1})$ in the usual way to bisect the length of 
$\gamma_n$. Let $\gamma_{A}$ and $\gamma_{B}$ be the restrictions 
of $\gamma_{n-1}$ to $[a_{n-1},c]$ and $[c,b_{n-1}]$ respectively. 
Suppose that $\gamma_{A} \geq k|\gamma_{A}(c)-\gamma_{A}(a_{n-1})|$ and 
$\gamma_{B} \geq k|\gamma_{B}(b_{n-1})-\gamma_{B}(c)|$. 
Then we have 
\begin{eqnarray*} 
|\gamma_{n-1}| & = & |\gamma_{A}|+|\gamma_{B}| \\ 
& \geq & k(|\gamma_{A}(c)-\gamma_{A}(a_{n-1})| + 
|\gamma_{B}(b_{n-1})-\gamma_{B}(c)|) \\ 
& \geq & 
k|\gamma_{A}(c)-\gamma_{A}(a_{n-1})+\gamma_{B}(b_{n-1})-\gamma_{B}(c)| \\ 
& = & k|\gamma_{n-1}(b_{n-1})-\gamma_{n-1}(a_{n-1})| 
\end{eqnarray*} 
which is a contradiction. 
Thus we must have either $|\gamma_{A}| < k|\gamma_{A}(c)-\gamma_{A}(a_{n-1})|$ 
or $|\gamma_{B}| < k|\gamma_{B}(b_{n-1})-\gamma_{B}(c)|$. 
We now set $\gamma_{n}$ to be either $\gamma_{A}$ or $\gamma_{B}$ accordingly. 

\QED 

We are now ready to prove the uniqueness of $\paths$-derivatives 
in the case where $\paths(X)$ is dense in $X$. 

\begin{theorem} 
Let $X \subseteq \Complex$ be compact and $\paths$ be a useful set of paths 
in $X$ such that $\paths(X)$ is dense in $X$. Then for $f \in C(X)$, any 
$\paths$-derivative of $f$ is unique. 
\end{theorem} 
\Proof 
Choose $f \in C(X)$. Suppose that $g, h \in C(X)$ are both 
$\paths$-derivatives of $f$. 
We have 
\[ 
\int_{\gamma} g(z) \wrt{z} = 
\int_{\gamma} h(z) \wrt{z} = 
f(z_{2})-f(z_{1}) 
\] 
for all $z_{1},z_{2} \in X$ and $\gamma \in \paths(z_{1},z_{2})$. 

Choose $z_{0} \in X$ and assume that $g(z_{0}) \neq h(z_{0})$. 
Now $g$ and $h$ are continuous, so there is an $R > 0$ and a $\delta > 0$ such 
that 
$|g(w)-h(w)| \geq \delta$ for all $w \in B_{X}(z_{0},R)$. 

Choose a path $\gamma \in \paths$ with $\Image{\gamma} \subseteq 
B_{X}(z_{0},R)$. 
We have seen that there is a sequence $(\gamma_{n})$ of sub-paths of $\gamma$ 
and a constant $k > 0$ such that $|\gamma_{n}| \rightarrow 0$ and 
\[|\gamma_{n}| < k|\gamma_{n}(b_{n})-\gamma_{n}(a_{n})|\] 
for each $n \in \Natural$. 
Now there must be a point $\alpha \in \cap_{n \in \Natural} 
\Image{\gamma_{n}}$ since 
the images are compact and nested. 
Clearly $g(\alpha)-h(\alpha) \neq 0$. For $z \in X$, write 
\[ 
g(z)-h(z) = g(\alpha)-h(\alpha)+r(z) 
\] 
where $r \in C(X)$ tends to 0 as $z \rightarrow \alpha$. 
Now choose $n \in \Natural$ such that $|r(z)| < \frac{\delta}{2k}$ for 
$z \in \Image{\gamma_{n}}$. 
We have 
\begin{eqnarray*} 
\int_{\gamma_{n}}(g(z)-h(z)) \wrt{z} 
& = & (g(\alpha)-h(\alpha))\int_{\gamma_{n}} \wrt{z} + \int_{\gamma_{n}} r(z) 
\wrt{z} \\ 
& = & (g(\alpha)-h(\alpha))(\gamma_{n}(b_{n})-\gamma_{n}(a_{n})) + 
\int_{\gamma_{n}} r(z) \wrt{z} 
\end{eqnarray*} 
Now 
\begin{eqnarray*} 
|(g(\alpha)-h(\alpha))(\gamma_{n}(b_{n})-\gamma_{n}(a_{n}))| 
& = & |g(\alpha)-h(\alpha)| \cdot |\gamma_{n}(b_{n})-\gamma_{n}(a_{n})| \\ 
& \geq & \delta \cdot |\gamma_{n}(b_{n})-\gamma_{n}(a_{n})| 
\end{eqnarray*} 
and 
\begin{eqnarray*} 
\left|\int_{\gamma_{n}} r(z) \wrt{z}\right| 
& \leq & |\gamma_{n}| \cdot \sup \{|r(z)| : z \in \Image{\gamma_{n}}\} \\ 
& \leq & k|\gamma_{n}(b_{n})-\gamma_{n}(a_{n})| \cdot \frac{\delta}{2k} \\ 
& = & |\gamma_{n}(b_{n})-\gamma_{n}(a_{n})| \cdot \frac{\delta}{2} 
\end{eqnarray*} 
Thus we have $\int_{\gamma_{n}}(g(z)-h(z)) \wrt{z} \neq 0$, but this is a 
contradiction. 
Hence the assumption that $g(z_{0}) \neq h(z_{0})$ must have been false, and 
so we have $g = h$. 

\QED 

\newcounter{unique1} 
\remembertheorem{unique1} 

Note that, even if $\paths(X)$ is not dense in $X$, it is clear that for any 
function 
$f \in \diff{\paths}{X}$, any two $\paths$-derivatives $g_{1}$ and $g_{2}$ of 
$f$ must agree on 
$\Closure{\paths(X)}$. This point will be crucial in the following development 
of the analytic properties of $\diff{\paths}{X}$. 

Note that the converse to the previous theorem is clear: 
if $\paths(X)$ is not dense in $X$ then every $f$ in 
$\diff{\paths}{X}$ has infinitely many $\paths$-derivatives. 
However, these $\paths$-derivatives will all agree on 
$\Closure{\paths (X)}.$ 

We now define the norm we need to make $\diff{\paths}{X}$ into a Banach 
algebra. 

\begin{definition} 
Let $X \subseteq \Complex$ be compact and $\paths$ be a non-empty, useful set 
of paths 
in $X$. For $f \in \diff{\paths}{X}$ we define 
\[\|f\| = \|f\|_{\infty} + \|g\|_{\Closure{\paths(X)}}\] 
where $g \in C(X)$ is any $\paths$-derivative of $f$ 
\end{definition} 

Note that $\|\cdot\|$ is well defined even when $\paths(X)$ is not dense in 
$X$ and $\paths$-derivatives are non-unique, because we know that 
any two $\paths$-derivatives of a function $f \in \diff{\paths}{X}$ 
\emph{do} agree on $\Closure{\paths(X)}$. 

\begin{theorem} 
Let $X \subseteq \Complex$ be compact and $\paths$ be a useful set of paths 
in $X$. Then $\diff{\paths}{X}$ is a normed space. 
\end{theorem} 
\Proof 
Clearly we have $\|f\| \geq 0$ for all $f \in \diff{\paths}{X}$, and 
$\|f\| = 0$ if and only if $f = 0$. 
Choose $f \in \diff{\paths}{X}$ and $\lambda \in \Complex$. 
Let $g \in C(X)$ be an $\paths$-derivative of $f$. 
We have already seen that $\lambda g$ is an $\paths$-derivative of $\lambda 
f$. 
We have: 
\begin{eqnarray*} 
\|\lambda f\| & = & \|\lambda f\|_{\infty} + \|\lambda 
g\|_{\Closure{\paths(X)}} \\ 
& = & |\lambda| \cdot \|f\|_{\infty} + |\lambda| \cdot 
\|g\|_{\Closure{\paths(X)}} \\ 
& = & |\lambda| \cdot \|f\| 
\end{eqnarray*} 
Choose $f_{1}, f_{2} \in \diff{\paths}{X}$, and $\paths$-derivatives $g_{1}$ 
and $g_{2}$ 
respectively. We have already seen that $g_{1}+g_{2}$ is an 
$\paths$-derivative of 
$f_{1}+f_{2}$. 
We have: 
\begin{eqnarray*} 
\|f_{1}+f_{2}\| 
& = & \|f_{1}+f_{2}\|_{\infty} + \|g_{1}+g_{2}\|_{\Closure{\paths(X)}} \\ 
& \leq & \|f_{1}\|_{\infty} + \|f_{2}\|_{\infty} + 
\|g_{1}\|_{\Closure{\paths(X)}} + \|g_{2}\|_{\Closure{\paths(X)}} \\ 
& = & \|f_{1}\| + \|f_{2}\| 
\end{eqnarray*} 

\QED 

We now show that $\diff{\paths}{X}$ is a Banach space. 

\begin{theorem} 
Let $X \subseteq \Complex$ be compact and $\paths$ be a useful set of paths in 
$X$. 
Then $\diff{\paths}{X}$ is complete. 
\end{theorem} 
\Proof 
Set $Y=\Closure{\paths(X)}$.
Let $(f_{n})$ be a Cauchy sequence in $\diff{\paths}{X}$. 
For each $n \in \Natural$, choose an $\paths$-derivative $g_{n}$ of $f_{n}$. 
Then $(f_{n})$ is Cauchy in $C(X)$ and 
$\left(g_{n}|_Y\right)$ is Cauchy in 
$C(Y)$, so these sequences converge uniformly.
Say $f_{n} \rightarrow f \in C(X)$ and $g_{n}|_Y \rightarrow g \in 
C(Y)$.

Extend $g$ to $\tilde{g} \in C(X)$ by the Tietze extension theorem. It is now
easy to check that $\tilde{g}$ 
is an $\paths$-derivative of $f$, so $f \in \diff{\paths}{X}$ and hence 
$\diff{\paths}{X}$ is complete. 
\QED

The last thing we have to do to show that $\diff{\paths}{X}$ is a Banach 
algebra is to show 
that $\|\cdot\|$ is an algebra norm. Fortunately this is not too difficult. 

\begin{theorem} 
Let $X \subseteq \Complex$ be compact and $\paths$ be a useful set of paths in 
$X$. 
Then the norm $\|f\| = \|f\|_{\infty} + \|g\|_{\Closure{\paths(X)}}$ (where 
$g$ is any 
$\paths$-derivative of $f$) is an algebra norm on $\diff{\paths}{X}$. 
\end{theorem} 
\Proof 
Choose $f_{1}, f_{2} \in \diff{\paths}{X}$, and $\paths$-derivatives $g_{1}$ 
and $g_{2}$ 
respectively. We have already seen that $g_{1}f_{2} + f_{1}g_{2}$ is an 
$\paths$-derivative 
of $f_{1}f_{2}$. 
We have 
\begin{eqnarray*} 
\|f_{1}f_{2}\| 
& = & \|f_{1}f_{2}\|_{\infty} + \|g_{1}f_{2} + 
f_{1}g_{2}\|_{\Closure{\paths(X)}} \\ 
& \leq & \|f_{1}\|_{\infty} \|f_{2}\|_{\infty} + 
\|g_{1}f_{2}\|_{\Closure{\paths(X)}} + \|f_{1}g_{2}\|_{\Closure{\paths(X)}} \\ 
& \leq & \|f_{1}\|_{\infty} \|f_{2}\|_{\infty} + 
\|g_{1}\|_{\Closure{\paths(X)}} \|f_{2}\|_{\infty} + \|f_{1}\|_{\infty} 
\|g_{2}\|_{\Closure{\paths(X)}} \\ 
& \leq & \|f_{1}\|_{\infty} \|f_{2}\|_{\infty} + 
\|g_{1}\|_{\Closure{\paths(X)}} \|f_{2}\|_{\infty} \\ 
& & + \|f_{1}\|_{\infty} \|g_{2}\|_{\Closure{\paths(X)}} + 
\|g_{1}\|_{\Closure{\paths(X)}} \|g_{2}\|_{\Closure{\paths(X)}} \\ 
& = & \left(\|f_{1}\|_{\infty} + 
\|g_{1}\|_{\Closure{\paths(X)}}\right)\left(\|f_{2}\|_{\infty} + 
\|g_{2}\|_{\Closure{\paths(X)}}\right) \\ 
& = & \|f_{1}\| \cdot \|f_{2}\| . 
\end{eqnarray*} 

\QED 

To avoid any complications arising from non-uniqueness of
$\paths$-derivatives, when we come to higher derivatives we will restrict 
attention to the case where $\paths (X)$ is dense in $X$. 

We next show that $\stddiff{X} \subseteq \diff{\paths}{X}$, and note 
conditions under which the inclusion is isometric. We also see the connection 
between the two kinds of derivative. 
\begin{theorem} 
Let $X \subseteq \Complex$ be compact and perfect and $\paths$ be a useful set 
of paths in $X$. 
Then $\stddiff{X} \subseteq \diff{\paths}{X}$. Indeed, 
for each $f \in \stddiff{X}$ the derivative (in the old sense) $f'$ is also an 
$\paths$-derivative of $f$. 
If $\Closure{\paths(X)} = X$ then the inclusion above is isometric. 
\end{theorem} 

\Proof 
Choose $f \in \stddiff{X}$. Then $f'$ (in the old sense) exists and is in 
$C(X)$. 
The Fundamental Theorem of Calculus for rectifiable paths gives us 
\[\int_{\gamma} f'(z) \wrt{z} = f(z_{2})-f(z_{1})\] 
for all $z_{1},z_{2} \in X$ and $\gamma \in \paths(z_{1},z_{2})$ and so $f'$ 
is an 
$\paths$-derivative of $f$. 
The rest is clear. 
\QED 

Note that when $\Closure{\paths(X)} = X$ the completion of 
$\stddiff{X}$ is simply its closure in $\diff{\paths}{X}$. 

We now introduce the new versions of the higher derivatives. As mentioned 
above, 
we will simplify matters by restricting attention to the case 
where 
$\Closure{\paths(X)} = X$. 
In view of the equality (in this setting) of the two kinds of derivatives when 
both are defined, 
we may safely use the notation $f'$ for the derivative of $f$ 
in either sense. 

Given such $X$ and $\paths$, it is clear how to define (inductively) 
the notion of {\it $n$-times $\paths$-differentiable} and the 
$n$th $\paths$-derivative of a function $f$. An easy induction 
using the above theorem shows that if $f$ is in $D^n(X)$ then 
$f$ is $n$-times $\paths$-differentiable and the old $n$th 
derivative $f^{(n)}$ is also the $n$th $\paths$-derivative of 
$f$. Thus we may use the notation $f^{(n)}$ for the new notion of derivative 
also. Moreover, in view of our earlier results, there 
is no problem (in this setting) in checking that the standard Leibniz formula 
is still valid for the new notion of $n$th derivative. 

We can now define spaces corresponding to the $D^n(X)$ spaces 
and the $D(X,M)$ spaces. We denote these new spaces by 
$\diffn{\paths}{n}{X}$ and ${\mathcal D}_\paths(X,M)$. 
For $f \in \diffn{\paths}{n}{X}$, we define 
\[ 
\|f\|_{n} = \sum_{k=0}^{n} \frac{\|f^{(k)}\|_{\Closure{\paths(X)}}}{k!} 
\] 
(with the usual convention that $f^{(0)} = f$). 
Similarly we define the norm on $\mathcal{D}_\paths(X,M)$ corresponding to 
that on $D(X,M)$. 

Because $\diff{\paths}{X}$ is complete, 
the new spaces are all Banach spaces and the old spaces are contained 
isometrically in the new spaces 
(the argument for this is the same as the proof that 
whenever $\stddiff{X}$ is complete then so are the $D^n(X)$ and $D(X,M)$ 
spaces). The spaces 
$\diffn{\paths}{n}{X}$ are always Banach algebras. When $M$ is 
an algebra sequence, $\mathcal{D}_\paths(X,M)$ is also a Banach algebra. The completions 
of the old spaces are simply their closures in the new spaces. 

In the next section we will investigate questions 
concerning the density or otherwise of $\stddiff{X}$ in 
$\diff{\paths}{X}$, along with some related questions of 
polynomial, rational and holomorphic approximation in these spaces and the higher 
derivative spaces.

\section{Approximation results} 

We will show that in many cases, $\diff{\paths}{X}$ 
is itself the completion of $\stddiff{X}$.
We begin with some cases where the two 
spaces are equal. In this first result, part of the conclusion (the fact that $\stddiff{X}$ is complete) 
was previously observed in
\cite{Honary+Mahyar1}. 
\begin{theorem} 
Let $X \subseteq \Complex$ be compact, perfect and the union of finitely many 
pointwise regular 
sets.
Let $L>0$, and let $\paths$ be a useful set of paths in $X$ which includes
all injective rectifiable paths with length $\leq L$ in $X$. Then 
$\stddiff{X} = \diff{\paths}{X}$ 
(and hence $\stddiff{X}$ is complete). 
\end{theorem} 
\Proof 
As we observed before, every rectifiable path is ``piecewise of length at most $L$'' 
and so we may assume that $\paths$ is, in fact, the set of {\it all} injective 
rectifiable paths in $X$.
By Lemma \recalltheorem{fin_pr} each component of $X$ is pointwise regular and 
so we have $\Closure{\paths(X)} = X$. 
Thus $\paths$-derivatives are unique and $\stddiff{X} \subseteq 
\diff{\paths}{X}$, 
the inclusion being isometric. 

Choose $f \in \diff{\paths}{X}$ and let $g \in C(X)$ be the 
$\paths$-derivative of $f$. 
Choose $a \in X$ and $(z_{n}) \subseteq X$ such that $z_n \neq a$ and $z_{n} \rightarrow a$. 
Then $z_{n}$ is eventually in the same component as $a$ (call the component 
$U$). So without loss of generality we can assume that $z_{n} \in U$ 
for every $n \in \Natural$. 

Since $U$ is pointwise regular, there is a $k(a)>0$ such that, for each $n \in \Natural$, there is an injective, rectifiable path $\gamma_{n}$ from $a$ to $z_{n}$ 
in $U$, 
such that $|\gamma_{n}| \leq k_{a}|a-z_{n}|$. 

Set 
\[ 
I_{n} = \frac{f(z_{n})-f(a)}{z_{n}-a} - g(a) 
= \left(\frac{1}{z_{n}-a}\int_{\gamma_{n}} g(z) \wrt{z}\right) - g(a) 
\] 
We have 
\[ 
g(a) = \frac{1}{z_{n}-a}\int_{\gamma_{n}} g(a) \wrt{z} 
\] 
for each $n \in \Natural$, so 
\begin{eqnarray*} 
|I_{n}| & = & \left| \frac{1}{z_{n}-a}\int_{\gamma_{n}} (g(z)-g(a)) \wrt{z} 
\right| \\ 
& \leq & \frac{|\gamma_{n}|}{|z_{n}-a|} \cdot \sup \{ |g(z)-g(a)| : z \in 
\Image{\gamma_{n}} \} \\ 
& \leq & k_{a} \cdot \sup \{ |g(z)-g(a)| : z \in \Image{\gamma_{n}} \} \rightarrow 0 
\end{eqnarray*} 
as $n \to \infty$. Thus 
\[ 
\lim_{n \rightarrow \infty} \frac{f(z_{n})-f(a)}{z_{n}-a} = g(a) 
\] 
as required. 

\QED 

\newcounter{pwreg_complete} 
\remembertheorem{pwreg_complete} 

Note that it is \emph{not} enough just to have $\paths$ being a useful set 
of paths in $X$ with $\Closure{\paths(X)} = X$, as the following example 
shows. 

\begin{example} 
Let $X$ be the unit square, $[0,1] \times [0,1] \subseteq \Complex$. 
Let $\paths$ be the set of paths of the form $\gamma : [a,b] \rightarrow X$, 
$\gamma(t) = k + ti$ for some set $[a,b] \subseteq [0,1]$. 
It is clear that $\paths$ is useful and $\Closure{\paths(X)} = X$, in fact 
$\paths(X) = X$. 

Define $f \in C(X)$ by $f(x+iy) = x$ for $x+iy \in X$. 
It is clear from the Cauchy-Riemann equations that $f \notin \stddiff{X}$. 
However $f \in \diff{\paths}{X}$, with (unique) $\paths$-derivative 
$g \in C(X)$ given by $g(x+iy) = 0$ for all $x+iy \in X$. 

In this example, $\stddiff{X}$ is a proper closed subalgebra of 
$\diff{\paths}{X}$. 
\end{example} 
\jremth{TooFewPaths}

In view of this example, it is worth investigating conditions on
$\paths$ which ensure that functions in $\diff{\paths}{X}$ are analytic on 
the interior of $X$. The following lemma and its immediate corollary give one 
class of useful sets of paths with this property.

\begin{lemma} 
Let $X \subseteq \Complex$ be compact and let
$\paths$ be the set of all 
injective rectifiable paths in $X$. Then every $f \in \diff{\paths}{X}$ 
is analytic on the interior of $X$. 
\end{lemma} 
\Proof 
Choose a point $z \in \mbox{int}(X)$. Then $B(z,r) \subseteq \mbox{int}(X)$ 
for some $r>0$. Set $Y=\Closure{B(z,r)}$. Let $\mathcal{G}$ be the set of paths in $\paths$ whose 
images are contained in $Y$. Clearly $\mathcal{G}$ is in fact 
the set of 
injective rectifiable paths in $Y$. 
Choose $f \in \diff{\paths}{X}$. 
Then $f|_Y \in \diff{\mathcal{G}}{Y}$. 
Since $Y$ is pointwise regular, by Theorem 
\recalltheorem{pwreg_complete} 
we have $\diff{\mathcal{G}}{Y} = \stddiff{Y}$. 
Thus $f|_Y \in \stddiff{Y}$
and hence $f$ is analytic at $z$. 

\QED 

\jremth{AnalyticOnInterior}

The following corollary is now immediate. 

\begin{corollary}
Let $X \subseteq \Complex$ be compact, Let $L>0$ and let
$\paths$ be a useful set of
paths in $X$ which includes all those injective rectifiable paths in $X$ which have length at most $L$. Then every 
$f \in \diff{\paths}{X}$ is analytic on the interior of $X$. 

\end{corollary}

\jremth{AnalyticOnInteriorCor}

%The following is a partial converse to Theorem %\recalltheorem{pwreg_complete}:

In order to show that the original spaces are often dense in the new spaces,
we will look at some related questions concerning polynomial and 
rational approximation. 
We first extend a definition from the original paper of Dales and Davie. 
(See \cite{Dales+Davie} and \cite{Honary} for some work 
on these spaces in the original setting.) 

\begin{definition} 
Let $X$ be a perfect, compact plane set 
Let $D$ be any of the normed algebras of functions on $X$ 
discussed in this paper such that 
$D$ includes all rational functions with poles of $X$. 
We define $D_R$ to be 
the closure in $D$ of the rational 
functions with poles off $X$ 
and $D_P$ to be the closure in $D$ of 
the polynomial functions. 
\end{definition} 

Curiously, it appears to be an open question whether or not $D_R$ 
is always equal to $D$ even restricting attention to 
the case where 
$D$ is $\stddiff{X}$. 
It is, however, easy to see that the continuous character space of 
$D_R$ is always equal to $X$ (recall that $D$ may be incomplete): the proof 
of Theorem 1.8 of \cite{Dales+Davie} goes through without need for 
any modifications. 
It is also elementary to see that whenever $\C \setminus X$ is connected then $D_P$= $D_R$: 
for example, this follows from the fact that the spectrum of the coordinate functional $Z$ must be the
same in the completions of $D_R$ and $D_P$.

\begin{theorem} 
Let $X$ be a compact plane set and suppose that $\paths$ 
is a useful set of paths in $X$ with 
$\Closure{\paths(X)} = X$. 
Let $D=\diff{\paths}{X}$. 
Consider the following subsets of $D$: 
$B$ is the set of all $f \in \diff{\paths}{X}$ such that 
the $\paths$-derivative of $f$ is the zero function; 
$A$ is the linear span of the idempotents in $D$ and 
$C$ is the closure in $D$ of $A$. 
Then $C$ is equal to the set of all functions 
in $C(X)$ which are constant on every component of $X$. 
Moreover we have 
$A\subseteq \stddiff{X}$, 
$C \subseteq B$ 
and $C \subseteq D_R$. 
\end{theorem} 
\Proof 
It is clear that all of the subsets of $D$ mentioned are in fact subalgebras 
of $D$ and 
that $A \subseteq B$ and $A \subseteq \stddiff{X}$. 
It is also clear that 
all of the idempotents in $C(X)$ are in $A$, and 
that $B$ is a closed subalgebra of $D$. 
Since the derivatives of all elements involved are $0$, 
$C$ is equal to the uniformly closed linear span in $C(X)$ 
of the idempotents in $C(X)$, and this is easily seen to be equal to 
the set of all functions 
in $C(X)$ which are constant on every component of $X$. 

Finally we turn to the rational approximation result. For this we need 
only prove that all of the idempotents are in $D_R$. 
This is, of course, immediate from the 
Shilov idempotent theorem, but can also be seen directly 
by extending any idempotent to be a function analytic on a neighbourhood 
of $X$ and applying Runge's theorem. 

\QED 

Note that Example \jrecth{TooFewPaths} shows that $C$ need not coincide with 
$B$ unless further conditions are placed on $\paths$. 

We now prove some closely related approximation results.

\begin{theorem} 
Let $X$ be a perfect, compact plane set 
such that $\C \setminus X$ is connected and $X$ has empty interior. Suppose that
$\paths$ is a useful 
set of paths in $X$ with the following properties: $\Closure{\paths(X)} = X$,
the set of lengths of paths in $\paths$ is bounded above 
and every pair of distinct 
points of $X$ which are in the same component of $X$ 
can be joined 
by a path in $\paths$.
Set $D=\diff{\paths}{X}$. 
Then $D=D_R=D_P$ and $\diff{\paths}{X}$ is the completion of $\stddiff{X}$. 
\end{theorem} 
\jremth{rational_approx}

\Proof 
It is clear that the second part of the conclusion follows from the first, and 
we have already mentioned that $D_P=D_R$ in this setting. 
We prove the rational approximation result. 
Set $L=\sup\{|\gamma| : \gamma \in \paths\}.$

Choose $f \in \diff{\paths}{X}$ with $\paths$-derivative $g \in C(X)$. 
By Mergelyan's theorem we can find a polynomial $p \in C(X)$ such that 
\[ 
\|p-g\|_{X} < \min\left\{\frac{\varepsilon}{3}, \frac{\varepsilon}{3 L}\right\} 
\] 
We have, for any path $\gamma_{0} \in \paths$, 
\begin{eqnarray*} 
\left|\int_{\gamma_{0}} (p(z)-g(z))\wrt{z}\right| 
& \leq & |\gamma_{0}| \cdot \|p-g\|_{X} \\ 
& < & \frac{\varepsilon|\gamma_{0}|}{3 L} 
\\ 
& \leq & \frac{\varepsilon}{3} 
\end{eqnarray*} 

Choose an analytic polynomial $F$ whose derivative is $p$.
Certainly we have $F \in D_R$. 

Since $F-f$ is uniformly continuous on $X$, we may choose 
$\delta > 0$ such that whenever $z,w \in X$ with 
$|z-w|< \delta$ then $|(F-f)(z)-(F-f)(w)| < \varepsilon/3$.

Noting that every component of $X$ is the intersection of the clopen sets containing it, 
by compactness we may choose components $K_1, K_2, \dots K_n$ of $X$ and pairwise disjoint 
clopen subsets $U_1, U_2, \dots, U_n$ of $X$ such that
for each $i$ we have 
$$U_i \subseteq \{z \in \Complex: \mbox{dist}(z,K_i)<\delta\},$$
and 
$$X=\cup_{i=1}^n U_i.$$

For each $i$, choose a point $z_i \in K_i$. Define
$h$ on $X$ as follows: for $z \in U_i$, $h(z)=F(z)- F(z_i) +f(z_i)$.
Then $h$ is $F$ plus 
a linear combination of idempotents, so $h \in D_R$.
We now look at $h-f$. First note that $h'=p$, so
$\|h'-g\|_{X} < \frac{\varepsilon}{3}$. We now wish to
estimate $\|h-f\|_X$. Although we have not assumed that 
$K_i \subseteq U_i$,
our choice of $\delta$ ensures that if $|F(z)- F(z_i) +f(z_i)-f(z)|$
is small on $K_i$ then $|h-f|$ is small on $U_i$.
Let $z \in K_i$. Choose any path $\gamma \in \paths$ from
$z_i$ to $z$. Then 
$$F(z)-F(z_i)+f(z_i)-f(z)=\int_\gamma{p(w)-g(w) \wrt w},$$
and so $|F(z)-F(z_i)+f(z_i)-f(z)| < \varepsilon/3$ for
$z \in K_i$. It follows from our choice of $\delta$
that $|h(z)-f(z)|< 2\varepsilon/3$ for 
$z \in U_i$ and so $|h-f|_X < 2\varepsilon/3$.
Thus $\|h-f\| < \varepsilon$.
\QED 

A similar theorem is valid when the interior is non-empty, provided that we ensure that the functions in $\diff{\paths}{X}$
are analytic on the interior of $X$. (This is of course also necessary for rational approximation). 
We saw some conditions which were sufficient for this above in Lemma \jrecth{AnalyticOnInterior} and 
Corollary \jrecth{AnalyticOnInteriorCor}. Here is one fairly general version of the result for 
$X$ with interior.

\begin{theorem} 
Let $X$ be a perfect, compact plane set 
such that $\C \setminus X$ is connected and let $r>0$.
Suppose that $\paths$ is a useful 
set of paths in $X$ 
with the following properties: $\Closure{\paths(X)} = X$, 
the set of lengths of paths in $\paths$ is bounded above, 
every pair of distinct 
points of $X$ which are in the same component of $X$ 
can be joined 
by a path in $\paths$ and $\paths$ includes all injective rectifiable paths in $X$ of length $\leq r$.
Set $D=\diff{\paths}{X}$. 
Then $D=D_R=D_P$ and $\diff{\paths}{X}$ is the completion of $\stddiff{X}$. 
\end{theorem} 
\jremth{rational_approx_int}

The proof is the same as that of \jrecth{rational_approx}
in view of the fact that, since $f$ is analytic on the 
interior of $X$, so is the $\paths$-derivative 
$g$ of $f$. As $g$ is continuous on $X$
we may still apply Mergelyan's theorem.

When $\C \setminus X$ is not connected then the polynomials can not be dense. If we attempt to 
imitate the above proofs using rational functions we hit the obstacle that it may not be 
possible to anti-differentiate
these. However, if a rational function is uniformly close to an $\paths$-derivative, we may obtain
good estimates on the residues at the poles and this may allow us to modify the rational function
slightly to obtain one which may be anti-differentiated. Here is one rational approximation result
valid for finitely connected $X$.

\begin{theorem}
Let $X$ be a perfect, compact plane set such that $\C \setminus X$ has only finitely many bounded
components,
say $U_1$, $U_2$, \dots, $U_n$. Choose one point $a_j$ in each of the bounded components $U_j$
$(1 \leq j \leq n).$
Suppose that $\paths$ is a useful set of paths in $X$ satisfying the conditions of Theorem
\jrecth{rational_approx_int} and, in addition, for each $j$ with $1 \leq j \leq n$ 
there is a closed curve $\gamma_j$ in $\paths$
with non-zero winding number about $a_j$.
Set $D=\diff{\paths}{X}$. 
Then $D=D_R$ and $\diff{\paths}{X}$ is the completion of $\stddiff{X}$. 
\end{theorem} 
\jremth{finitely-connected}

\Proof
The proof is again similar to that of Theorem \jrecth{rational_approx}.  
Let $f \in \diff{\paths}{X}$ and let $g$ be the $\paths$-derivative
of $f$. Then $g$ is continuous on $X$ and analytic on the interior of $X$, so it is standard
(see for example \cite{Gamelin})
that $g$ may be uniformly approximated on $X$ by a sequence of rational functions, say
$r_k$.
By Runge's theorem we may further assume that the poles of the rational functions $r_k$
all lie in
$\{a_1, a_2, \dots,a_n\}$.
As $\int_{\gamma_j} g(z)~\rd z = 0$ for each $j$, the residue at each $a_j$ of $r_k$
tends to $0$ as $k \to \infty$. Thus we may modify the sequence $r_k$ (subtracting
rational functions with simple poles in $\{a_1, a_2, \dots,a_n\}$ if necessary)
to show that 
$g$ may be uniformly approximated on $X$ by anti-differentiable rational functions.
The remainder of the proof is identical to that of Theorem \jrecth{rational_approx}, using such a
rational function $r$ in place of the polynomial $p$ used there and taking $F$ to be a rational 
anti-derivative of $r$.
\QED

These theorems cover many cases where completeness has previously
been
an issue, for example the simple sets and the radially self absorbing
sets
considered in \cite{FLO} (we define these below), and the combs
and stars considered in \cite{Bland_Thesis}.

Recall the following definition from \cite{FLO}.

\begin{definition}
Let $X \subseteq \Complex$ be non-empty and compact. Then $X$ is
\emph{radially
self-absorbing} if, for every $r > 0$ we have $X \subseteq \mbox{int}(rX)$.
\end{definition}

We conclude this section by transferring to our new spaces a result 
about holomorphic
approximation for radially self absorbing sets, 
originally proved for the $D(X,M)$ spaces in
\cite{FLO} (Lemma 3.1). The bulk of 
the proof is identical, but we need to use Lemma \jrecth{AnalyticOnInterior}. 

\newcounter{diffanalytic} 
\remembertheorem{diffanalytic} 

\begin{theorem} 
Let $X \subseteq \Complex$ be compact and radially self-absorbing. 
Let $\paths$ be the set of injective rectifiable paths in $X$. 
Let $M$ be a sequence of positive numbers. 
Set 
\[ 
S = \{f \in D(X,M) : f \mbox { extends to be analytic on
a neighbourhood of } X\}. 
\] 
Then $S$ is dense in $\mathcal{D}_{\paths}(X,M)$. 
\end{theorem} 
\Proof 
Note that $\Closure{\paths(X)} = X$ and so $D(X,M)$ embeds isometrically in 
$\mathcal{D}_{\paths}(X,M)$.

Choose $f \in \mathcal{D}_{\paths}(X,M)$. Then by Lemma
\jrecth{AnalyticOnInterior}, 
$f$ is analytic on $\mbox{int}(X)$. 
For $n \in \Natural$ and $z \in \Complex$, set 
$g_{n}(z) = \frac{nz}{n+1}$
and set
$F_{n} = f \circ g_{n}$. 
Then $F_{n}$ is analytic on $\left(\frac{n+1}{n}\right)\mbox{int}(X)$,
which is a neighbourhood of $X$, and so
$F_{n}|_{X} \in D^{\infty}(X)$. 

Set $f_{n} = F_{n}|_{X}$. 
We have, for all $n \in \Natural$ and all $k \geq 0$, 
$\|f_{n}^{(k)}\|_\infty 
\leq 
\|f^{(k)}\|_\infty$.
Thus $f_{n} \in D(X,M)$ for each $n \in \Natural$. 
Now $\|f_{n}^{(k)} - f^{(k)}\|_\infty \rightarrow 0$ as $n 
\rightarrow \infty$. 
Hence, by dominated convergence for series, $\|f_{n}-f\| \rightarrow 0$ 
as $n \rightarrow \infty$. 

\QED 

{\bf Remarks}
This holomorphic approximation result also shows that the new space is the 
completion of the old. The same result is, of course, valid for the $\diffn{\paths}{n}{X}$ spaces 
(with a slightly easier proof). 

In the algebra setting, the completeness of our $\mathcal{D}_{\paths}(X,M)$ algebras 
allows us to apply the holomorphic functional calculus to the coordinate 
functional $Z$, as in Corollary 3.2 of \cite{FLO}, 
to see that the polynomials are dense in our spaces:
the same proof goes through without the need for any changes. 
Of course, this result from \cite{FLO} follows from our result,
without need for the completeness assumption made there. This 
eliminates the need to appeal to the Runge argument given in 
Section 5 of that paper to cover the possibility that the normed
algebras concerned
might be incomplete.

\section{Open Problems} 

We conclude with some open problems. 
%Fairly easy to see that pointwise regular implies locally connected.

1. Does there exist a compact plane set $X$ and a non-analytic sequence $M$ 
such that the rational functions with poles off $X$ are not dense in 
$D(X,M)$? 

We have mentioned a small number of positive results on polynomial and 
rational approximation, but in general this problem is wide open. In 
particular the answer for $D(X,M)$ is apparently not known for the 
`square annulus' obtained by 
deleting an open square from the middle of a compact square. 
(Note that in view of these open problems, some authors have worked directly 
with the closures in these spaces of the set of rational functions instead.) 

2. If $X$ is a radially self-absorbing set, is $\stddiff{X}$ already complete? 
More generally, suppose that $X$ is the closure of a bounded, connected open subset of 
$\Complex$. Is $\stddiff{X}$ already complete? 

3. Let $X$ be a compact plane set and let $\paths$ be the set of all injective, rectifiable paths
in $X$. Suppose that $\Closure{\paths(X)}=X$. Is it always true that
$\diff{\paths}{X}$ is the completion of $\stddiff{X}$? Is it always true that the rational functions
are dense in $\diff{\paths}{X}$?

Note that here the answer for the square annulus is easily seen to be yes, 
using Theorem \jrecth{finitely-connected}. More generally, any function 
$f \in \diff{\paths}{X}$
which may be extended to have continuous first-order partial derivatives on 
a neighbourhood of $X$ and whose $\overline\partial$ derivative vanishes on
$X$ may be approximated in $\diff{\paths}{X}$ by rational functions
(see \cite{Verdera}). However, even functions in $\stddiff{X}$ need not in
general have such extensions. Our question on rational approximation
is equivalent to the question 
of whether a dense set of functions may be so extended.

\vskip 0.3cm

%William Bland
% Macrovision

% 2830, De La Cruz Blvd

% Santa Clara

% CA 95050

% email: William.Bland@abstractnonsense.com

%and

%Joel F. Feinstein

$^*$ {\bf Corresponding author}

{\sf  School of Mathematical Sciences

 University of Nottingham

 Nottingham NG7 2RD, England

 email: Joel.Feinstein@nottingham.ac.uk
\vskip 0.3cm

2000 Mathematics Subject Classification: 46J15, 46E10}


\begin{thebibliography}{99} 

\bibitem{Ali_Thesis} Abdullah K. S. Ali, 
``Some properties of Banach Function Algebras'', 
Ph.D. Thesis, 
University of Nottingham, 1995. 

\bibitem{Apostol} Tom M. Apostol, 
``Mathematical Analysis'', 
Addison-Wesley, 1974. 

\bibitem{Behrouzi} F. Behrouzi,
``Homomorphisms of certain Banach function algebras'',
Proc. Indian Acad. Sci. Math. Sci. 112 (2002), 331-336.

\bibitem{Bland_Thesis} William J. Bland, 
``Banach Function Algebras And Their Properties'', 
Ph.D. Thesis, 
University of Nottingham, 2001. 

\bibitem{Conway} John B. Conway, 
``Functions of a complex variable. I.'', 
Springer, New York, 1978. 
\bibitem{Garthsbook} H. G. Dales, 
``Banach algebras and automatic continuity'', 
LMS Monographs 24, 
Clarenden Press, Oxford, 2000. 
\bibitem{Dales+Davie} H. G. Dales and A. M. Davie, 
``Quasianalytic Banach function algebras'', 
J. Funct. Anal. {\bf 13} (1973), 
28-50. 

\bibitem{Dales+McClure} H.G. Dales and J.P. McClure, ``Completion of normed 
algebras of 
polynomials'', 
J. Austral. Math. Soc. (A) {\bf 20} (1975), 504-510. 

\bibitem{FLO} J. F. Feinstein, H. Lande and A. G. O'Farrell, 
``Approximation and extension in normed spaces of infinitely differentiable 
functions'', 
J. London Math. Soc. (2) {\bf 54} (1996), 
541 - 556. 
\bibitem{Feinstein+Kamowitz} Joel F. Feinstein and Herbert Kamowitz, 
``Endomorphisms of Banach algebras of infinitely differentiable 
functions on compact plane sets'', 
J. Funct. Anal. {\bf 173} (2000), 61 - 73. 
\bibitem{Gamelin}
T.W. Gamelin, ``Uniform algebras'', Prentice-Hall, New Jersey, 1969.

\bibitem{Honary} 
T.G. Honary, 
``Relations between Banach function algebras and their uniform closures'', 
Proc. Amer. Math. Soc. 109 (1990), no. 2, 337--342. 

\bibitem{Honary+Mahyar1} 
T.G. Honary, H. Mahyar, 
``Approximation in Lipschitz algebras of infinitely differentiable 
functions'', 
Bull. Korean Math. Soc. 36 (1999), no. 4, 629--636. 
\bibitem{Honary+Mahyar2} 
T.G. Honary, H. Mahyar, 
``Approximation in Lipschitz algebras'', 
Quaest. Math. 23 (2000), no. 1, 13--19 


\bibitem{Jarosz} 
K. Jarosz, 
``${\rm Lip}\sb {\rm Hol}(X,\alpha)$'', 
Proc. Amer. Math. Soc. {\bf 125} (1997), 3129--3130. 

\bibitem{Kamowitz} H. Kamowitz, 
``Endomorphisms of Banach algebras of infinitely differentiable functions'', 
Banach Algebras '97 (Blaubeuren), 
De Gruyter, Berlin, (1998), 
273-285. 

\bibitem{OFar} A.G. O'Farrell, ``Polynomial approximation of smooth 
functions'', 
J. London Math. Soc. (2) 28 (1983), 496-506 

\bibitem{Rudin} Walter Rudin, 
``Real and complex analysis'', 
McGraw-Hill, 1987. 

\bibitem{Verdera} J. Verdera, ``On $C^m$ rational approximation'',
Proc. Amer. Math. Soc. {\bf 97} (1986), 621--625.
\end{thebibliography}
\end{document}